\documentstyle{article}
\newtheorem{prop}{Proposition}[section]
\newtheorem{th}[prop]{Theorem}

\newtheorem{lm}[prop]{Lemma}
\newtheorem{df}[prop]{Definition}
\newtheorem{ex}[prop]{Example}
\newtheorem{remark}[prop]{Remark}
\newtheorem{question}[prop]{Question}
\newtheorem{fact}[prop]{Fact}

\newcommand{\Ker}{{\rm{Ker}}}

\newcommand{\bsquare}{\hbox{\rule{6pt}{6pt}}}
\newcommand{\proof}[1]{\noindent{\bf Proof}\hspace{0.3cm}{#1}\hfill\bsquare 
\vspace{0.5cm}\par}

\date{October 23, 2017}

\begin{document}

\title{Zariski K3 surfaces}
\author{Toshiyuki Katsura\thanks{Partially supported by JSPS Grant-in-Aid 
for Scientific Research (B) No. 15H03614} 
~and Matthias Sch\"utt\thanks{Partially supported by ERC StG 279723.}
}

\maketitle

\begin{abstract}
We construct Zariski K3 surfaces of Artin invariant 1, 2 and 3 in many characteristics.
In particular, we prove that any supersingular Kummer surface is Zariski if $p\not\equiv 1$ mod $12$.
Our methods combine different approaches such as quotients by the group scheme $\alpha_p$,
Kummer surfaces, and automorphisms of hyperelliptic curves.
\end{abstract}

\section{Introduction}

Let $k$ be an algebraically closed field of characteristic $p$, and $X$ be
an algebraic variety of dimension $n$ over $k$. $X$ is said to be unirational if
there exists a dominant rational map from ${\bf P}^n$ to $X$.
Unirational algebraic curves are automatically rational,
and it is a classical fact that the same holds for complex algebraic surfaces
as a consequence of Castelnuovo's criterion for rationality.
In positive characteristic, however, this is no longer true
as was first shown by Zariski \cite{Zariski}.
In essence, this is due to the impact of inseparable maps
which leads to the following definition:

\begin{df}
A (non-rational) algebraic surface $S$ is called a Zariski surface if there exists a 
purely inseparable dominant rational map ${\bf P}^{2}\to S$
of degree $p$.
\end{df}

In this sense, Zariski surfaces can be considered the minimal non-rational unirational surfaces.
Note that automatically ${\rm H}^2_{et}(S,{\bf Q}_\ell)$ is spanned by algebraic cycles,
i.e.~Zariski surfaces are supersingular.
This leads to the question to what extent the converse may hold true
(as initiated by Shioda \cite{Shioda-uni}).
Here we concentrate on supersingular K3 surfaces;
partly this is due to their striking history in this problem,
dating back to the first discoveries of Zariski and unirational K3 surfaces by Artin \cite{A} and Shioda \cite{Shioda-p=3},
but mostly because  the unirationality of supersingular K3 surfaces is known by now 
outside characteristic $3$
by
Rudakov--Shafarevich
\cite{RS}  and Liedtke \cite{Liedtke}.

\begin{question}
\begin{enumerate}
\item[(i)]
Is any supersingular K3 surface a Zariski surface? 
\item[(ii)]
In particular,
is any supersingular Kummer surface  a Zariski surface?
\end{enumerate}
\end{question}

While the first question has an affirmative answer
if the characteristic $p = 2$ (cf. Rudakov-Shafarevich
\cite[p.~151]{RS} where the proof of the corollary indeed implies the surfaces to be Zariski),
a general answer might be too much to ask for at this time.
Meanwhile the second question was prompted by the initial result of Shioda \cite{Shioda-some}
that all supersingular Kummer surfaces are unirational.
We will develop an affirmative answer to \textit{(ii)} for 75\% of all characteristics:

\begin{th}
\label{th}
Let $p>2$ such that $p\not\equiv 1$ mod $12$.
Then any supersingular Kummer surface in characteristic $p$ is a Zariski surface.
\end{th}

In addition, we will provide a plentitude of new Zariski K3 surfaces.
Our results are summarized in the following theorem
where only part \textit{(i)} seems to have been known before
(see, e.g., Katsura \cite{K1}).

\begin{th}
\label{thm}
There are Zariski K3 surfaces of Artin invariant $\sigma$
over an algebraically closed field of characteristic $p$ under the following conditions:
\begin{enumerate}
\item[(i)]
$\sigma=1$ and $p\not\equiv 1  \;{ \rm mod } \; 12$.
\item[(ii)]
$\sigma=2$ and $p\not\equiv 1, 49\;  { \rm mod }\;  60$.
\item[(iii)]
$\sigma=3$ and $p\equiv 3,5\; { \rm mod } \;7$.
\end{enumerate}
\end{th}

The proofs of Theorems \ref{th} and \ref{thm} proceed by explicit geometric constructions.
We combine different approaches such as quotients by the group scheme $\alpha_p$,
Kummer surfaces, lattice theory, and automorphisms of elliptic and hyperelliptic curves.
We supplement the theorems with additional results in two directions:
an isolated Zariski  K3 surface of Artin invariant $\sigma=1$ over ${\mathbf F}_{13}$
(Example \ref{ex:13} -- this surface also is Kummer),
and  an abundance of Zariski elliptic surfaces (Lemma \ref{lem:ell}, Remark \ref{rem:ell}).

\section{Preliminaries on supersingular abelian surfaces}
\label{s:prel}

For later use, we start by reviewing  parts of the theory of abelian surfaces in positive characteristic.
Throughout the paper, we fix an algebraically closed field $k$ of characteristic $p > 0$.
An abelian surface $A$ is said to be \emph{supersingular} (resp. \emph{superspecial}) if it is isogenous 
(resp. isomorphic) to
a product of two supersingular elliptic curves.
In terms of the N\'eron-Severi lattice,
this phrases as 
\[
\rho(A) = \mbox{rank NS}(A) = 6
\;\;\; \mbox{ with discriminant } -p^{2\sigma}.
\]
Here $\sigma$ is called \emph{Artin invariant} and equals $1$ if $A$ is superspecial, and $2$ otherwise.
By definition a superspecial abelian surface is supersingular (cf. Oort \cite{Oo}),
and a superspecial abelian surface is unique up to isomorphism (cf. Shioda \cite{S}). 
In this section, we recall some results on the N\'eron-Severi 
group of the superspecial abelian surface. 

Let $E$ be a supersingular elliptic curve defined over $k$.
We consider the superspecial abelian surface $E_{1} \times E_{2}$
with $E_{1} = E_{2} = E$. We  denote by $O_{E}$
the zero point of $E$. We put $X = E_{1}\times\{O_{E_2}\} + 
\{O_{E_1}\} \times E_{2}$, which is a principal polarization on $E_{1} \times E_{2}$.
By abuse of notation, we sometimes denote the fibers $E_{1}\times\{O_{E_2}\}$  
(resp.  $\{O_{E_1}\} \times E_{2}$) by $E_{1}$ (resp. by $E_{2}$).
We set ${\cal O} = {\rm End}(E)$ and $B = {\rm End}^{0}(E)=
{\rm End}(E)\otimes {\bf Q}$.
$B$ is a quaternion division algebra over the rational number field ${\bf Q}$
with discriminant $p$, and ${\cal O}$ is a maximal order of $B$ (cf. Mumford \cite{M}, Section 22). 
For an element $a \in B$, we denote by $\bar{a}$ the image under the canonical involution.
We have a natural identification of ${\rm End}(E_{1} \times E_{2})$ 
with the ring ${\rm M_{2}}({\cal O})$
of two-by-two matrices with coefficients in ${\cal O}$.
Here, the action of $\left(
\begin{array}{cc}
\alpha  & \beta \\
\gamma & \delta
\end{array}
\right) \in {\rm M_{2}}({\cal O})$ is given by
$$
\begin{array}{rccc}
\left(
\begin{array}{cc}
\alpha  & \beta \\
\gamma & \delta
\end{array}
\right) : &  E_{1} \times E_{2} & \longrightarrow & E_{1} \times E_{2} \\
    & (x, y)   &\mapsto & (\alpha (x) + \beta (y), \gamma (x) + \delta (y)).
\end{array}
$$

By a divisor $L$ we usually mean the divisor class
represented by $L$ in ${\rm NS}(E_{1} \times E_{2})$ if confusion is unlikely to occur.
With this convention, a divisor $L$ yields  a homomorphism
$$
\begin{array}{cccc}
   \varphi_{L} :  &　E_{1} \times E_{2}  &\longrightarrow  & {\rm Pic}^{0}(E_{1} \times E_{2}) \\
      & x & \mapsto & T_{x}^{*}L - L,
\end{array}
$$
where $T_{x}$ is the translation by $x \in E_{1} \times E_{2}$ (cf Mumford \cite{M}).
We set 
$$
     H = \left\{
\left(
\begin{array}{cc}
\alpha  & \beta \\
\gamma & \delta
\end{array}
\right) \in {\rm M}_{2}({\cal O})
~\mid
~\alpha, \delta \in {\bf Z},~\gamma, \beta \in {\cal O},~\gamma = \bar{\beta}
\right\}.
$$
Note that for an automorphism $g$ of $E_{1} \times E_{2}$, we can regard $g$ as an element of
${\rm M}_{2}({\cal O})$. By this identification, we have ${}^{t}\bar{g}$ as an element of
${\rm M}_{2}({\cal O})$.
We will use the following theorem which is well-known to specialists
 (see Katsura \cite{K2}, for instance).

\begin{th}\label{intersection}
The homomorphism
$$
\begin{array}{cccc}
j : & {\rm NS}(E_{1} \times E_{2}) & \longrightarrow  & H \\
   & L  &\mapsto  & \varphi_{X}^{-1}\circ\varphi_{L}
\end{array}
$$
is bijective.  By this correspondence, we have
$$
j(E_{1}) = \left(
\begin{array}{cc}
0  & 0 \\
0 & 1
\end{array}
\right),~
j(E_{2}) = \left(
\begin{array}{cc}
1  & 0 \\
0 &  0
\end{array}
\right).
$$
For $L_{1}, L_{2} \in {\rm NS}(E_{1} \times E_{2})$ such that
$$
  j(L_{1}) =
\left(
\begin{array}{cc}
\alpha_{1}  & \beta_{1} \\
\gamma_{1} & \delta_{1}
\end{array}
\right),~
  j(L_{2}) =
\left(
\begin{array}{cc}
\alpha_{2}  & \beta_{2} \\
\gamma_{2} & \delta_{2}
\end{array}
\right),
$$
the intersection number  $L_{1}\cdot L_{2}$ is given by
$$
L_{1}\cdot L_{2} = \alpha_{2}\delta_{1} + \alpha_{1}\delta_{2} - \gamma_{1}\beta_{2}
-\gamma_{2}\beta_{1}.
$$
In particular, for $L\in {\rm NS}(E_{1} \times E_{2})$ such that 
$j(L) = \left(
\begin{array}{cc}
\alpha  & \beta \\
\gamma & \delta
\end{array}
\right)
$
we have
$$
L^{2} = 2\det \left(
\begin{array}{cc}
\alpha  & \beta \\
\gamma & \delta
\end{array}
\right), \;\;\;
L\cdot E_{1}  = \alpha,
\;\;\;
L\cdot E_{2} = \delta .
$$
We have also  $j(nD)=nj(D)$ for an integer $n$.

For $L_{1}, L_{2} \in {\rm NS}(E_{1} \times E_{2})$ with $j(L_{1}) = g_{1}$ and 
$j(L_{2}) = g_{2}$ and for an automorphism $g$ of $E_{1} \times E_{2}$, we have
 $g^{*}L_{1} \equiv L_{2}$ if and only if ${}^{t}\bar{g} g_{1}g = g_{2}$.
\end{th}

Let $E$ be a supersingular elliptic curve defined over ${\bf F}_p$. Such an elliptic curve
exists for any $p > 0$ (cf. Waterhouse \cite{W}). We denote by $F$ 
the relative Frobenius morphism of $E$.
For the local-local group scheme $\alpha_p$
of rank $p$, we have ${\rm End} (\alpha_p) \cong k$. Therefore, for $(i, j) \in k^2$
we have an inclusion 
$$
\epsilon : \alpha_p \stackrel{(i, j)}{\longrightarrow} \alpha_p \times \alpha_p \subset E \times E.
$$
We assume $i/j \notin {\bf F}_{p^{2}}, j\neq 0$. Then, by Oort \cite{Oo}, 
the quotient surface 
$$A = (E \times E)/\epsilon (\alpha_p)
$$ 
is not superspecial.
Let 
$$
\pi : E \times E \longrightarrow A
$$ 
be the quotient map. Considering the dual abelian surface
$A^t$ of $A$ and the dual homomorphism $\pi^t$, we have a commutative diagram
for $D\in {\rm NS}(A)$:
$$
\begin{array}{ccccc}
     E \times E &  \stackrel{\varphi_{\pi^{*}D}}{\longrightarrow} & E \times E & 
\stackrel{\varphi_X^{-1}}{\longrightarrow} & E \times E\\
        \pi \downarrow    &       &\uparrow \pi^t  &  &\\
          A   &  \stackrel{\varphi_D}{\longrightarrow}  & A^t. &  & 
\end{array}
$$
Here, we identify the dual abelian surface $(E \times E)^t$ with $E \times E$.
We set
$$
H' =
\left\{
\left(
\begin{array}{cc}
\alpha  & \beta \\
\gamma & \delta
\end{array}
\right) \in {\rm M}_{2}({\cal O})
~\mid
\begin{array}{l}
~\alpha, \delta \in {\bf Z},~\gamma, \beta \in {\cal O},~\gamma = \bar{\beta}, \\
p\mid \alpha,~
p\mid \delta,~ \beta \in F{\cal O} = F{\rm End (E)}
\end{array}
\right\} 
\subset H \subset {\rm M}_{2}({\cal O}).
$$
\begin{prop}\label{H'}
The homomorphism
$$
   j\circ \pi^* : {\rm NS}(A) \longrightarrow H
$$
induces an isomorphism $ j\circ \pi^* : {\rm NS}(A) \longrightarrow H'$
of additive groups.
\end{prop}

\proof{The proof is essentially the same as the one in Ibukiyama-Katsura-Oort \cite{IKO}, 
Proposition 2.4.1.
Since we have $\dim_k{\rm Hom}(\alpha_p, A) = \dim_k{\rm Hom}(\alpha_p, A^t) = 1$,
the subgroup schemes which are isomorphic to $\alpha_p$ are unique in $A$ and $A^t$,
respectively. Therefore, we have $\varphi_D^{-1}(\alpha_p) \supset \alpha_{p}$ 
for $D \in {\rm NS}(A)$. This implies that
$\Ker \varphi_{\pi^{*}D} \supset \alpha_p \times \alpha_p$. Setting
$$
j({\pi^{*}D}) = \left(
\begin{array}{cc}
\alpha  & \beta \\
\gamma & \delta
\end{array}
\right),
$$
and considering that $\alpha$ and $\delta$ are integers, we have $p\mid \alpha$,
$p\mid \delta$ and $\beta, \gamma  \in F{\cal O}$, that is, $j({\pi^{*}D}) \in H'$.

Conversely, for $j(G) \in H'$ with $G \in {\rm NS}(E\times E)$, we have
$\Ker {\varphi_{G}} \supset \alpha_p \times \alpha_p$. Using the notation in Mumford \cite{M},
for any subgroup scheme $\epsilon(\alpha_p)$ of $\alpha_p \times \alpha_p$
we have $e^D(\epsilon(\alpha_p), \epsilon(\alpha_p)) = 0$. Therefore, using the descent theory
in Mumford \cite{M}, there exists a divisor $G' \in {\rm NS}(A)$ such that $\pi^*(G') = G$.
Hence $ j\circ \pi^* : {\rm NS}(A) \longrightarrow H'$ is surjective.
}

Let $\beta_i$ ($i = 1, 2, 3, 4$) be a basis of ${\cal O}$ over ${\bf Z}$. Then,
$$
\left\langle \left(
\begin{array}{cc}
0  & 0 \\
0 & 1
\end{array}
\right),~
\left(
\begin{array}{cc}
1  & 0 \\
0 & 0
\end{array}
\right),~
\left(
\begin{array}{cc}
0  & \beta_i\\
\bar{\beta}_i & 0
\end{array}
\right)~  (i = 1, 2, 3, 4)
\right\rangle
$$
is a basis of $H$. There exist $D_i$ ($i = 1, 2, 3, 4$) such that 
$j(D_i) = \left(
\begin{array}{cc}
0  & \beta_i\\
\bar{\beta}_i & 0
\end{array}
\right)$.
By Theorem \ref{intersection}, $\langle E_1, E_2, D_1, D_2, D_3, D_4\rangle$ is a basis of ${\rm NS}(E\times E)$.
Therefore, the determinant of the Gram matrix $M$ of the given basis is equal to $- p^2$.
The following proposition is well-known (cf. Ogus \cite{O}). We give an elementary proof
for it.
\begin{prop}\label{Abelian-Artin}
The Artin invariant of $A$ is equal to 2.
\end{prop}
\proof{A basis of $H'$ is given by
$$
\left\langle \left(
\begin{array}{cc}
0  & 0 \\
0 & p
\end{array}
\right),~
\left(
\begin{array}{cc}
p  & 0 \\
0 & 0
\end{array}
\right),~
\left(
\begin{array}{cc}
0  & F\beta_i\\
\overline{F\beta}_i & 0
\end{array}
\right)~  (i = 1, 2, 3, 4)\right\rangle
$$
By Proposition \ref{H'}, there exist divisors $D'_i$ $(i = 1, 2, \cdots , 6)$ 
in ${\rm NS}(A)$ such that
$$
\begin{array}{l}
j(\pi^*(D'_{1}) ) =\left(
\begin{array}{cc}
0  & 0 \\
0 & p
\end{array}
\right),~
j(\pi^*(D'_{2}) ) =\left(
\begin{array}{cc}
p  & 0 \\
0 & 0
\end{array}
\right),\\
j(\pi^*(D'_{i +2}) )  = \left(
\begin{array}{cc}
0  & F\beta_i\\
\overline{F\beta}_i & 0
\end{array}
\right)~  (i = 1, 2, 3, 4).
\end{array}
$$
Here, $\pi^*(D'_{i}) = pE_i$ ($i = 1, 2$).
Since the Gram matrix $((\pi^*D'_i, \pi^*D'_j)) =( p(D'_i, D'_j))$, we have
$$
\begin{array}{l}
p^6\det ((D'_i, D'_j)) = \det ((\pi^*D'_i, \pi^*D'_j)) \\
= 
\det \left(
\begin{array}{cc}
\begin{array}{cc}
      (pE_1)^2 & (pE_1, pE_2) \\
    (pE_2, pE_1) &  (pE_1)^2 
\end{array} 
    & 0   \\

0  &  ((\pi^*(D'_{i +2}), \pi^*(D'_{j +2}))_{1\leq i, j \leq 4}
\end{array}
\right)\\
= 
\det \left(
\begin{array}{cc}
\begin{array}{cc}
      0 & p^2 \\
    p^2 &  0
\end{array} 
    & 0   \\

0  &  (p(D_{i}, D_{j}))_{1\leq i, j \leq 4}
\end{array}
\right) = p^8\det M = - p^{10}.
\end{array}
$$
Therefore, we have $\det ((D'_i, D'_j)) = - p^4$, that is, the Artin invariant of $A$
is equal to 2.
}

\section{Generalized K3 surfaces}
\label{s:gen}

It is standard (outside characteristic $2$) to associate to an abelian surface
a K3 surface by means of the Kummer quotient.
In this section, we shall discuss different constructions
which have the benefit of two compatibilities,
both with the constructions in Section \ref{s:prel} and with purely inseparable base change
as required for Zariski surfaces.
Throughout we retain the notation from Section \ref{s:prel}.

First, we assume 
$p \equiv 2 ~({\rm mod}~3)$
and consider the supersingular elliptic curve $E$ with j-invariant zero defined by
\begin{eqnarray}
\label{eq:j=0}
E: \;\;\;    y^2 + y= x^3.
\end{eqnarray}
$E$ is endowed with an automorphism $\tau$ of order 3 defined by
$$
\tau: \;\;\;      x \mapsto \omega x,\;\; y\mapsto y
$$
where $\omega$ denotes a primitive cube root of unity. 

\begin{remark}
\label{rem:3}
Some readers may be more familiar with the Weierstrass form
\[
y^2 = x^3 -1
\]
for $E$. Outside characteristic $2$, both models are isomorphic,
but (\ref{eq:j=0}) comes with the advantage of being valid in characteristic $2$ as well
while often also yielding simpler equations.
\end{remark}

\begin{fact}
The endomorphism ring of $E$ can be represented as
\begin{eqnarray}
\label{eq:OO}
{\rm End}(E) =
{\cal O} =
 {\bf Z} \oplus {\bf Z}F 
\oplus {\bf Z}\tau \oplus {\bf Z}(1 + F)(2 + \tau)/3.
\end{eqnarray}
\end{fact}

\proof{
In characteristic $p>2, p\equiv 2$ mod $3$, this is Katsura \cite{K1}, Lemma 5.4
(building on Ibukiyama \cite{Ibu}).
Meanwhile, for $p=2$ one may verify directly that the ring ${\cal O}$ from (\ref{eq:OO}) is isomorphic
to the maximal order in the quaternion algebra $\left(\frac{-1,-1}{\mathbf Q}\right)$ of discriminant $2$,
for instance using that both elements $\omega_4, \omega_3-\omega_4$
(in the notation of loc. cit.) have square $-1$. This clearly implies the claim.
}

In End$(E)$, we have the relations $F\tau = \tau^2F$, $\bar{\tau} = \tau^2$, $\bar{F} = -F$ and 
$$\overline{(1 + F)(2 + \tau)}/3 = 1 - \{(1 + F)(2 + \tau)/3\}.
$$
For the sake of simplicity, we set $\eta = (1 + F)(2 + \tau)/3$.
The multiplication is given by the following table:

\begin{center}
\small{
\begin{tabular}{|c|c|c|c|c|}
\hline
     & $1$ & $F$ & $\tau$ & $\eta$ \\
\hline
$1$ & $1$ & $F $ & $\tau$ & $\eta$ \\
\hline
$F$ & $ F$ & $-p$ & $3\eta -2F -\tau -2$ & $\eta - 2(p + 1)/3 - (1 + p)\tau/3$ \\
\hline
$\tau$ & $\tau$ & $-3\eta + F + \tau + 2$ & $-\tau -1$ & $-2\eta + F + \tau + 1$\\
\hline
$\eta$ & $\eta$ &$(2 - p)/3 + F + (p + 1)\tau/3 -\eta$  & $\eta - F -1$  & $\eta -(p + 1)/3$  \\
\hline
\end{tabular}
}
\end{center}

The automorphism $\tau \times \tau$ acts on $E \times E$. Since 
$\tau \times \tau$ preserves the subgroup scheme  $\epsilon ({\alpha_p})$,
$\tau \times \tau$ induces an automorphism $\theta$ on the quotient 
\[
A=(E\times E)/\alpha_p.
\]
$\tau \times \tau$ has 9 isolated fixed point on $E\times E$
and $\pi$ is a finite purely inseparable morphism, we see that $\theta$
has also 9 isolated fixed points on $A$. By a local calculation,
the action of $\theta$ at the fixed points is given 
by $(s, t) \mapsto (\omega s, \omega^2 t)$. Therefore, the singularities
of the quotient space $A/\langle \theta\rangle$ are all rational double points
of type $A_2$. Considering the action of $\theta$ on the vector space 
${\rm H}^{0}(A, \Omega_A^1)$, we see that $\theta$ acts
symplectically on
a non-zero regular 2-form. We denote by ${\rm GKm} (A)$ the nonsingular complete 
minimal model of $A/\langle \theta\rangle$.
Since $A$ is a supersingular abelian surface, 
we conclude that ${\rm GKm}(A)$ is a supersingular K3 surface.

\begin{th}\label{theorem3}
Assume $p \equiv 2 ~({\rm mod}~3)$.
Under the notation as above,  the Artin invariant 
of ${\rm GKm}(A)$ is equal to 2.
\end{th}

\begin{remark}
A similar result in characteristic 2 is
given in Schr\"oer \cite{S}.
\end{remark}

\proof{We consider the quotient morphism $\pi : E \times E \longrightarrow A$.
Since we have a commutative diagram
$$
\begin{array}{rccc}
      \tau \times \tau : & E \times E & \longrightarrow & E \times E \\
           &   \downarrow &         & \downarrow \\
         \theta : & A   & \longrightarrow & A,
\end{array}
$$
it should be clear that
for $D\in {\rm NS}(A)$, 
$\theta^*(D) = D$ if and only if $(\tau \times \tau)^*(\pi^*D) = \pi^*D$.
We calculate the invariants $(H')^{\langle \tau \times \tau\rangle}$. By Theorem \ref{intersection},
the action of $\tau \times \tau$ on 
$\left(
\begin{array}{cc}
\alpha  & \beta \\
\gamma & \delta
\end{array}
\right) \in H'$ is given by
$$
\left(
\begin{array}{cc}
\alpha  & \beta \\
\gamma & \delta
\end{array}
\right) \mapsto 
\left(
\begin{array}{cc}
\bar{\tau}  & 0 \\
0 & \bar{\tau}
\end{array}
\right)
\left(
\begin{array}{cc}
\alpha  & \beta \\
\gamma & \delta
\end{array}
\right)
\left(
\begin{array}{cc}
\tau  & 0 \\
0 & \tau
\end{array}
\right).
$$
$\left(
\begin{array}{cc}
0  & 0 \\
0 & p
\end{array}
\right)$ and 
$\left(
\begin{array}{cc}
p  & 0 \\
0 & 0
\end{array}
\right)$ 
are invariant under this action. Now, let $F(a + b F + c\tau + d\eta)$
be an element of $F{\cal O}$ with $a, b , c, d \in {\bf Z}$ which satisfies 
$\tau^2(F(a + b F + c\tau + d\eta) \tau = F(a + b F + c\tau + d\eta)$.
Since $\tau^2 F = F\tau$, 
we have $\tau (a + b F + c\tau + d\eta) \tau = (a + b F + c\tau + d\eta)$.
Then, using the multiplication table, we have
$$
a + c + d = 0,~ 2a -c + d = 0, \mbox{that is}, c = a/2,~d = -3a/2.
$$
Hence, a basis of the invariant space $(H')^{\langle \tau \times \tau\rangle}$ is given by
$$
\begin{array}{c}
\left(
\begin{array}{cc}
0  & 0 \\
0 & p
\end{array}
\right),~
\left(
\begin{array}{cc}
p  & 0 \\
0 & 0
\end{array}
\right),~
\left(
\begin{array}{cc}
0  & p \\
p & 0
\end{array}
\right), \\
\left(
\begin{array}{cc}
0  & F(2 + \tau - 3\eta) \\
\overline{F(2 + \tau - 3\eta)} & 0
\end{array}
\right) = 
\left(
\begin{array}{cc}
0  & 2p + p\tau \\
2p + p\bar{\tau} & 0
\end{array}
\right).
\end{array}
$$
We take divisors $G_i$ ($i = 1, 2, 3, 4$) on $E\times E$ such that $j(G_i)$ correspond 
to elements of  this basis in this order. Then, there exist divisors $G'_i$ 
 ($i = 1, 2, 3, 4$) on $A$ such that $\pi^*(G'_i) = G_i$.  
The divisors $G'_i$  ($i = 1, 2, 3, 4$) form a basis of 
the invariant space ${\rm NS}(A)^{\langle \theta\rangle}$. By Theorem \ref{intersection},
the Gram matrix of the $G_i$ ($i = 1, 2, 3, 4$) is given by
$$
\left(
\begin{array}{cccc}
0  & p^2  &0 & 0\\
p^2  & 0 & 0 &  0\\
0  & 0 & -2p^2  &  -3p^2  \\
 0 & 0 &  -3p^2  &  -6p^2
\end{array}
\right).
$$
Its determinant is equal to $-3p^8$. Therefore, the determinant
of the Gram matrix of the basis $G'_i$  ($i = 1, 2, 3, 4$) is equal to $-3p^4$.
Hence, in a similar way to Katsura \cite{K1}, Lemma 5.8, we conclude
that the discriminant of ${\rm NS}({\rm GKm}(A))$ is equal to $-p^4$ as claimed.
}

Now, we assume $p \equiv 3 ~({\rm mod}~4)$. 
We consider a supersingular elliptic curve $E$ with j-invariant $1728$ defined by
\begin{eqnarray}
\label{eq:1728}
E:\;\;\;   y^2 = x^3 -x.
\end{eqnarray}
$E$ possesses an automorphism $\tau$ of order 4 defined by
$$
\tau: \;\;      x \mapsto - x,\;\; y\mapsto iy.
$$
Here $i$ is a primitive fourth root of unity. 
In this case, we have
$$
{\cal O} ={\rm End}(E) = {\bf Z} \oplus {\bf Z}\tau 
\oplus {\bf Z}(1 + F)/2 \oplus {\bf Z}\tau (1 + F)/2
$$
(cf. Katsura \cite{K1}, Lemma 5.3).
One checks that $\tau^2 = -1$, $F\tau = \tau^3F$, $\bar{\tau} = \tau^3$, 
$\bar{F} = -F$.
For the sake of simplicity, we set $\eta = \tau (1 + F)/2$.
The multiplication is given by the following table:

\begin{center}
\small{
\begin{tabular}{|c|c|c|c|c|}
\hline
     & $1$ & $\tau$ & $(1 + F)/2$ & $\eta$ \\
\hline
$1$ & $1$ & $\tau $ & $(1 + F)/2$ & $\eta$ \\
\hline
$\tau$ & $ \tau$ & $-1$ & $\eta$ & $- \eta$ \\
\hline
$(1 + F)/2$ & $(1 + F)/2$ & $\tau - \eta$ & $-(1 + p)/4 + (1 + F)/2$ & $(1 + p)\tau/4$\\
\hline
$\eta$ & $\eta$ &$-1 + (1 + F)/2$  & $-(1 + p)\tau/4 + \eta$  & $-(1 + p)/4$  \\
\hline
\end{tabular}
}
\end{center}

The automorphism $\tau \times \tau$ acts on $E \times E$. Since 
$\tau \times \tau$ preserves the subgroup scheme  $\epsilon ({\alpha_p})$,
$\tau \times \tau$ induces an automorphism $\theta$ on $A$.
By a similar method as above (also see Katsura \cite{K1}, Lemma 5.8),
we conclude that the nonsingular complete minimal model ${\rm GKm}(A)$ of 
$A/\langle \theta\rangle$ is a supersingular K3 surface.

A basis of the invariant space $(H')^{\langle \tau \times \tau\rangle}$ is given by
$$
\begin{array}{c}
\left(
\begin{array}{cc}
0  & 0 \\
0 & p
\end{array}
\right),~
\left(
\begin{array}{cc}
p  & 0 \\
0 & 0
\end{array}
\right),~
\left(
\begin{array}{cc}
0  & p \\
p & 0
\end{array}
\right), 
\left(
\begin{array}{cc}
0  & p\tau \\
- p\tau  & 0
\end{array}
\right).
\end{array}
$$
We take divisors $G_i$ ($i = 1, 2, 3, 4$) on $E\times E$ such that $j(G_i)$ correspond 
to elements of  this basis in this order. Then, there exist divisors $G'_i$ 
 ($i = 1, 2, 3, 4$) on $A$ such that $\pi^*(G'_i) = G_i$.  
The divisors $G'_i$  ($i = 1, 2, 3, 4$) form a basis of 
the invariant space ${\rm NS}(A)^{\langle \theta\rangle}$. By Theorem \ref{intersection},
the Gram matrix of $G_i$ ($i = 1, 2, 3, 4$) is given by
$$
\left(
\begin{array}{cccc}
0  & p^2  &0 & 0\\
p^2  & 0 & 0 &  0\\
0  & 0 & -2p^2  &  0  \\
 0 & 0 &  0  &  -2p^2
\end{array}
\right).
$$
The determinant of this matrix equals $-4p^8$. It follows that the determinant
of the Gram matrix of the basis $G'_i$  ($i = 1, 2, 3, 4$) is equal to $-4p^4$.

By a similar method to the above and Katsura \cite{K1}, Lemma 5.8,
we obtain the following theorem. We omit the details.

\begin{th}\label{theorem4}
Assume $p \equiv 3 ~({\rm mod}~4)$. 
Under the notation as above, 
the Artin invariant of ${\rm GKm}(A)$ is equal to 2.
\end{th}

\section{Zariski surfaces}

We now turn to the problem of Zariski surfaces,
in particular for K3 surfaces.
To this end, we let $C$ be a non-singular complete model of the algebraic curve defined by
\begin{equation}\label{equation1}    
C: \;\;\; ~y^2 + y = x^\ell
\end{equation}
with an integer  $\ell \geq 3$. Let $\zeta$ be a primitive $\ell$-th root of
unity (so $p\not | \ell$). Then, $C$ has an automorphism of order $\ell$ defined by
\begin{equation}\label{equation2}
\quad \quad \quad  \tau : \;\; x \mapsto \zeta x,\; y \mapsto y.
\end{equation}

\begin{remark}
As in the case of $\ell=3$ from Remark \ref{rem:3}, a more standard equation may consist in
\[
y^2 = x^\ell - 1.
\]
\end{remark}

We put $Y = C_1 \times C$ with $C_1 = C$. 
We assume that the defining equation for $C_1$ is given 
by 
$$
C_1: \;\;\; y_{1}^2 + y_1 = x_1^\ell.
$$
Then, $\tau \times \tau$ acts on $Y$ as in (\ref{equation2}).

\begin{lm}\label{rational}
$Y/\langle \tau \times \tau \rangle$ is a rational surface.
\end{lm}
\proof{
The group $G = \langle \tau\times \tau \rangle$ acts on
the function field of $Y$ via its natural action on $k(x_1, y_1, x, y)$. We set
$$    
     z = x/x_{1}.
$$
Then, the invariant field $k(Y)^G$ is given by
$k(y, y_{1}, z)$ with the relation 
$$
z^\ell(y_1^2 + y_1) =y^2 + y.
$$
This endows $Y/G$ with the structure of a conic fibration over ${\bf P}^1$ with the parameter $z$
and section $(0,0)$, say.
Therefore, this is a rational surface.
}

\begin{prop}\label{Zariski}
Assume $p \equiv i ~\mbox{mod}~ \ell$ $(2 \leq i \leq p -1)$. 
Then, $Y/\langle \tau \times \tau^i \rangle$ is a Zariski surface.
\end{prop}
\proof{We have the following commutative diagram:
$$
\begin{array}{rrl}  
            Y & \stackrel{\tau \times \tau}{\longrightarrow} & Y \\
  {\rm id} \times F   \downarrow &      & \downarrow {\rm id} \times F \\
          Y  & \stackrel{\tau \times \tau^i}{\longrightarrow} & Y.
\end{array}
$$
Therefore, we have a purely inseparable rational map 
$Y/\langle \tau \times \tau \rangle \longrightarrow Y/\langle \tau \times \tau^i \rangle$.
Since $Y/\langle \tau \times \tau \rangle$ is a rational surface by Lemma \ref{rational}, we conclude that 
$Y/\langle \tau \times \tau^i \rangle$ is a Zariski surface.
}

Now let $X$ be a K3 surface defined over $k$. $X$ is said to be \emph{supersingular}
if the Picard number satisfies
$$
\rho (X) = b_{2}(X) = 22
$$ 
with $b_{2}$ the second Betti number.
If $X$ is supersingular, then the discriminant of the N\'eron-Severi group 
is of the form $-p^{2\sigma}$ with a positive integer $\sigma\in\{1,\dots,10\}$
called the \emph{Artin invariant} 
(cf Artin \cite{A}). 
In analogy with abelian surfaces, a supersingular K3 surface with $\sigma = 1$ 
is called superspecial. If $p \neq 2$, it is known that any superspecial K3 surface
is isomorphic to the Kummer surface ${\rm Km}(E\times E)$ 
with $E$ a supersingular elliptic curve, and that any supersingular K3 surface
with Artin invariant 2 is isomorphic to a Kummer surface 
${\rm Km}((E\times E)/\alpha_{p})$ 
with $E$ a supersingular elliptic curve and a suitable embedding 
$\alpha_{p} \rightarrow E\times E$ (cf. Oort \cite{Oo}, Ogus \cite{O}, Theorem 7.10,
and see also Shioda \cite{Shi}). 
Here, $\alpha_{p}$ is the local-local
group scheme of rank 1 discussed in Section \ref{s:prel}.

We give a natural proof for the following known theorem (see Katsura \cite{K1}, for instance)
which will later be generalized in different directions.

\begin{th}\label{1}
Assume $p \not\equiv 1~({\rm mod}~12)$.
Then, the supersingular K3 surface with Artin invariant $\sigma=1$ is a Zariski surface.
\end{th}

Note that Theorem \ref{1} covers \textit{(i)} of Theorem \ref{thm}.
For the first missing case of characteristic $p=13$, see Example \ref{ex:13}.

\medskip

\proof{First, we assume $p \equiv 2$ mod  $3$ and take $E$ from (\ref{eq:j=0})
-- or from (\ref{equation1}) with automorphism
$\tau$ as in (\ref{equation2}) with $\ell = 3$.
Then, the minimal resolution of 
$(E\times E)/\langle \tau \times \tau^2 \rangle$ is isomorphic to
a Kummer surface ${\rm Km} (E\times E)$ (cf. Katsura \cite{K1}).
Hence, by Proposition \ref{Zariski}, ${\rm Km} (E\times E)$ is 
a Zariski surface with $C = E$ and $i = 2$.

Secondly, we assume $p \equiv 3~({\rm mod}~4)$. 
It follows that $E$ from (\ref{eq:1728})
 is a supersingular elliptic curve which is
isomorphic to the elliptic curve defined by (\ref{equation1})
with an automorphism $\sigma$ as in (\ref{equation2})  with $\ell = 4$.
Then, the minimal resolution of 
$(E\times E)/\langle \sigma \times \sigma^3 \rangle$ is isomorphic to
a Kummer surface ${\rm Km} (E\times E)$ (cf. Katsura \cite{K1});
hence it has Artin invariant $\sigma=1$.
As before, it is also a Zariski surface by Proposition \ref{Zariski} (with $i=3$).
This completes our proof.
}

In fact, we have already enough information to prove a big portion
(in terms of the congruence $p\not\equiv 1,49$ mod $60$) of Theorem \ref{thm} \textit{(ii)}.

\begin{th}\label{2}
Assume $p \not\equiv 1~({\rm mod}~12)$.
Then, there exists a Zariski supersingular K3 surface with Artin invariant 2.
\end{th}

\proof{
Assume $p\equiv 2~({\rm mod}~3)$ and let $E$ denote the supersingular elliptic curve from (\ref{eq:j=0}), (\ref{equation1})
with automorphism
$\tau$ as in (\ref{equation2}) with $\ell = 3$.
Consider the quotient morphism $\pi: E \times E \longrightarrow A$ as in Section 3.
This morphism induces 
$$
(E \times E)/\langle \tau \times \tau \rangle  \longrightarrow A/\langle \theta\rangle.
$$
The minimal model of $A/\langle \theta\rangle$ is ${\rm GKm}(A)$.
By Proposition \ref{rational}, $(E \times E)/\langle \tau \times \tau \rangle$ is
rational, 
so GKm$(A)$ is Zariski,
while by Theorem \ref{theorem3}, the Artin invariant of ${\rm GKm}(A)$ is equal to 2.

Secondly, assume $p\equiv 3~({\rm mod}~4)$. Fix the supersingular elliptic curve
$E$ from (\ref{equation1})
with  automorphism
$\tau$ as in (\ref{equation2}) with $\ell = 4$.
Recall that $E$ is  isomorphic to the elliptic curve
from (\ref{eq:1728}) with automorphism $\tau$ introduced in Section 3.
Arguing as above, but with  Theorem \ref{theorem4} instead,
we complete our proof.
%
%
}

\section{Kummer surfaces with Artin invariant 2}

In the previous section, we showed that there exist Zariski supersingular K3 surfaces
with Artin invariant 2 if $p \not\equiv 1~({\rm mod}~12)$. In this section,
we will prove Theorem \ref{th} by
showing directly that any supersingular K3 surface with Artin invariant 2
is Zariski if $p \not\equiv 1~({\rm mod}~12)$. 
Since Rudakov and Shafarevich already showed
that all supersingular K3 surfaces are Zariski in characteristic 2 
(Rudakov-Shafarevich \cite{RS}), we may assume $p \neq 2$.
We start with a few preparations.


\begin{lm}\label{action}
Let $A$ be an abelian surface, and let $\tau$ be an automorphism of $A$
of order $m$ ($m > 2$, $m\neq 4$). Let $\zeta$ be a primitive $m$-th root of unity.
Assume that $\tau$ acts as the multiplication by $\zeta$ 
on the vector space ${\rm H}^0(A, \Omega_{A}^{1})$ 
of regular 1-forms on $A$.
Then, the quotient surface $A/\langle \tau\rangle$ is rational.
\end{lm}


\proof{Let $X$ be a nonsingular model of $A/\langle \tau\rangle$.
Then, there exists a dominant rational map 
$$\varphi : A \longrightarrow X.
$$
Suppose that the Albanese variety ${\rm Alb}(X)$ of $X$ is nontrivial. 
Then, there exists a non-zero regular 1-form on ${\rm Alb}(X)$. Pulling back 
the regular 1-form to $X$, we have a non-zero regular 1-form $\omega$ on $X$.
Then, $\varphi^*(\omega)$ is a non-zero $\tau$-invariant regular 1-form on $A$,
which contradicts our assumption. Therefore, we see the irregularity $q(X) = 0$.

Now, suppose there exists a non-zero regular 2-ple 2-form $\Omega$ on $X$. Then,
$\varphi^*\Omega$ gives a non-zero $\tau$-invariant regular 2-ple 2-form on $A$.
However, since 
$$
{\rm H}^0(A, (\Omega^2_{A})^{\otimes 2}) \cong (\wedge^2{\rm H^0}(A, \Omega_{A}^{1}))^{\otimes 2}\cong k,
$$ 
the action of $\tau$ on the space of regular 2-ple 2-forms on $A$ is given by
 multiplication by $\zeta^4$, which is not $1$ by assumption, a contradiction. Hence, we conclude that
$X$ is rational by Castelnuovo's criterion of rationality as in Zariski \cite{Zariski}.
}

\begin{remark}
It is posible to weaken the assumption of Lemma \ref{action} as long as none of the induced actions of $\tau$
on differential forms is trivial.
In order to cover the case $m=4$, however, we will need to throw in some extra work in \ref{ss}.
\end{remark}

Let $X$ be an algebraic surface with $\dim {\rm H}^2(X, {\cal O}_X) = 1$. 
We assume that the formal Brauer group of $X$ is prorepresentable by
a one-dimensional formal group (cf. Artin-Mazur \cite{AM}). We denote the formal Brauer group by $\Phi_X$.

\begin{lm}\label{height}
Let $X$ (resp. $Y$) be an algebraic surface with $\dim {\rm H}^2(X, {\cal O}_X) = 1$ (resp. $\dim {\rm H}^2(Y, {\cal O}_Y)= 1$). 
Assume that their formal Brauer groups are prorepresentable by one-dimensional formal groups $\Phi_X$, $\Phi_Y$, respectively.
Moreover, assume there exists a dominant separable rational map $f: Y \longrightarrow X$ 
such that the degree of $f$ is prime to $p$.
Then, the height of $\Phi_X$ is equal to the height of $\Phi_Y$.
\end{lm}

\proof{
Since the formal Brauer group is stable under blowing-up by Artin-Mazur \cite{AM}, 
we have a homomorphism $f^{*}: \Phi_{X} \longrightarrow \Phi_{Y}$.
We also have a non-zero homomorphism 
$f^{*}: {\rm H}^2(X, {\cal O}_X) \longrightarrow 
{\rm H}^2(Y, {\cal O}_{Y})$, which is an isomorphism 
since each space is 1-dimensional and $f$ is separable whose degree is prime to $p$. 
Since ${\rm H}^2(X, {\cal O}_{X})$ 
(resp. ${\rm H}^2(Y, {\cal O}_Y)$) is the tangent space 
of $\Phi_{X}$ (resp. $\Phi_{Y}$), the homomorphism 
$f^{*}: \Phi_{X} \longrightarrow \Phi_{Y}$ is nontrivial.
Therefore, the height of $\Phi_{X}$ is equal to the height of $\Phi_{Y}$.
}

Now, we recall the theory of \emph{a-number} (for the details of a-number for algebraic varieties,
 see van der Geer-Katsura \cite{GK}, Definition 2.1).
 For a nonsingular complete algebraic surface $X$,
we denote by ${\rm H}_{dR}^2(X)$ the second De Rham cohomology group of $X$.
From here on, we consider only algebraic surfaces such that 
the Hodge-to-De Rham spectral sequence is degenerate at $E_1$-term.
(For instance, this holds  if the characteristic satisfies $p > 2$ and if $X$ can be lifted to the Witt ring $W_2(k)$,
and in particular for K3 surfaces.)
Then, we have the Hodge filtration
$$
{\rm H}_{dR}^2(X) = {\rm F}_0 \supset {\rm F}_1 \supset {\rm F}_2 \supset 0
$$
such that $F_0/F_1 = {\rm H}^2(X, {\cal O}_X)$, $F_1/F_2 = {\rm H}^1(X, \Omega^1_X)$
and $F_2 = {\rm H}^0(X, \Omega^2_X)$.
The absolute Frobenius map $F$ acts on ${\rm H}_{dR}^2(X)$, and 
the kernel of $F$ is ${\rm F}_1$. Therefore, we have an injective map
$$
   F : {\rm H}^2(X, {\cal O}_X) \longrightarrow {\rm H}_{dR}^2(X).
$$
Then, the a-number $a(X)$ of $X$ is defined by
$$
   a(X) = \max\{i \mid ({\rm Im} F) \cap {\rm F}_i \neq 0\}.
$$
Here, $\max$ means the largest number in the set.
Note that if $X$ is an abelian surface, $a(X)$ coincides with the usual a-number
defined by Oort \cite{Oo} (cf. Katsura \cite{GK}, Proposition 2.2).
For instance, a supersingular abelian surface $A$ of Artin invariant $2$ has a-number $1$,
since $\alpha_p$ embeds uniquely into $A$.

\begin{lm}\label{a-number}
Let $X$ (resp. $Y$) be an algebraic surface with $\dim {\rm H}^2(X, {\cal O}_X) = 1$ (resp. $\dim {\rm H}^2(Y, {\cal O}_Y)= 1$). 
Assume there exists a dominant separable rational map $f: Y \longrightarrow X$ 
such that the degree of $f$ is prime to $p$.
Then, we have $a(X) = a(Y)$.
\end{lm}
\proof{If neccesary, we blow up $Y$, and we may assume that $f$ is a morphism.
We have a commutative diagram:
$$
\begin{array}{ccc}
   {\rm H}^2(Y, {\cal O}_X) & \stackrel{F}{\longrightarrow} &{\rm H}_{dR}^2(Y)\\
    f^* \uparrow   &       &  f^* \uparrow   \\
   {\rm H}^2(X, {\cal O}_X) & \stackrel{F}{\longrightarrow} & {\rm H}_{dR}^2(X).
\end{array}
$$
Here, the second up-arrow preserves the Hodge filtrations. 
Since the degree of $f$ is prime to p, the first up-arrow is an isomorphism.
Therefore, we have $a(X) = a(Y)$.
}

\begin{lm}\label{a-number of K3}
Let $X$ be a supersingular K3 surface with Artin invariant 2
in characteristic $p>2$. 
Then the a-number of $X$
is equal to 1.
\end{lm}
\proof{
Since $X$ is a supersingular K3 surface with Artin invariant 2, 
$X$ is isomorphic to a Kummer surface ${\rm Km}(B)$, where $B$ is a supersingular
abelian surface with Artin invariant 2 (as we have mentioned before). Clearly, the a-number of $B$ is equal to 1.
Since we have a dominant separable rational map $f: B \longrightarrow X$ of degree 2 (which is prime to $p$),
Lemma \ref{a-number} implies that $a(X) = 1$.
}

\subsection{Proof of Theorem \ref{th} for $p\equiv 2$ mod $3$}
\label{ss3}

First we assume $p \equiv 2~({\rm mod}~3)$ and $p \neq 2$.
Let $E$ be the supersingular elliptic curve defined by
$y^2 + y = x^3$ and $\tau$ the automorphism defined by $x \mapsto \omega x$,
$y \mapsto y$ with $\omega$ a primitive cube root of unity. We take an immersion 
$$
\epsilon : \alpha_p \stackrel{(i, j)}{\longrightarrow} \alpha_p \times \alpha_p \longrightarrow  E \times E
$$
with $i/j \not\in {\bf F}_{p^2}$. Then, as we already showed, the automorphism 
$\tau\times \tau$ of $E\times E$ induces an order 3 automorphism $\theta$ 
of $A = (E \times E)/\epsilon (\alpha_p)$ and the nonsingular minimal model ${\rm GKm}(A)$
of $A/\langle \theta \rangle$ is a supersingular K3 surface with Artin invariant 2.

In spirit, our approach follows closely Shioda \cite[proof of Thm. 4.2, 4.3]{Shi}.
In fact, it can be adapted for the congruence class $p\equiv 3$ mod $4$ (in \ref{ss})
and goes roughly as follows:
\begin{enumerate}
\item
set up a smooth covering corresponding to the quotient map $A\to{\rm GKm}(A)$;
\item
translate the information into lattices
which thus carry over to any supersingular K3 surface $X$ of   Artin invariant 2;
\item
recover a cover of $X$ which leads to a supersingular abelian surface $Y$ of   Artin invariant 2;
\item
facilitate the unique embedding $\alpha_p\hookrightarrow Y$ to derive that $X$ is  a Zariski surface.
\end{enumerate}
The numbering of the subsections to follow reflects the above steps.

\subsubsection{}
First, we recall how to construct ${\rm GKm}(A)$, following Katsura \cite{K1}, Section 5.
Since $\theta$ has 9 fixed points on $A$, we blow up at these 9 points: 
$$\psi_1 : A_1 \longrightarrow A.
$$
 We denote by $G_i$ ($i = 1, \ldots, 9$)
the exceptional curves.
Then, $\theta$ induces an automorphism $\theta_1$ on $A_1$
which has 2 fixed points on each exceptional curve. We once more blow up
these 18 fixed points: 
$$
\psi_2 : A_2 \longrightarrow A_1.
$$ 
Abusing notation, we denote again by $G_i$ the proper transform of $G_i$, and
by $D_i$, $F_i$  ($i = 1, 2, \ldots, 9$) the  exceptional curves 
such that 
$$(D_i, G_i) = (F_i, G_i) = 1,\;\; (D_i,F_i)=0 \;\;\; (i = 1, 2, \ldots, 9).
$$
We have $D_i^2 = -1$, $F_i^2 = -1$, and $G_i^2 = -3$.
The automorphism $\theta_1$ on $A_1$ induces an automorphism
$\theta_2$ on order 3 on $A_2$
which acts as identity map on $D_i$ and $F_i$,
and induces an automorphism of order 3 on $G_i$.
Therefore, the quotient surface $A_2/\langle \theta_2\rangle$ is nonsingular
and we have a diagram:
$$
\begin{array}{ccccl}
     A_2 & \stackrel{\psi_2}{\longrightarrow} & A_1 & \stackrel{\psi_1}{\longrightarrow} & A\\
     \pi \downarrow & & & & \\
     A_2/\langle \theta_2\rangle & \stackrel{\psi_3}{\longrightarrow} & {\rm GKm}(A). & & 
\end{array}
$$
Here, $\pi$ is the quotient map and $\psi_3$ will be described momentarily.
We set $D'_i = \pi(D_i)$, $F'_i = \pi(F_i)$ and $G'_i = \pi(G_i)$. Then, we have
$\pi^*(D'_i) = 3D_i$, $\pi^*(F'_i) = 3F_i$ and $\pi^*(G'_i) = G_i$. Using these relations,
we have $(D'_i)^2 = -3$, $(F'_i)^2 = -3$ and $(G'_i)^2 = -1$.
The morphism $\psi_3$ then simply is the blowing-down of the 9 curves $G'_i$ ($i = 1, \ldots, 9$).

\subsubsection{}
Since $\pi : A_2 \longrightarrow  A_2/\langle \theta_2\rangle$ is a cyclic covering 
of degree 3 with smooth ramification locus comprising the $D_i$ and $F_i$, 
it induces an effective branch divisor 
$$
R = \sum_{i=1}^9 (D_i'+2F_i')
$$ 
which is divisible by 3
in ${\rm NS}(A_2/\langle \theta_2 \rangle)$; that is, there exists a divisor $R'$ 
on $A_2/\langle \theta_2 \rangle$ such that $3R' = R$.
(The coefficients of the components of $R$ guarantee that $R'$ has integral intersection number with each $G_i'$.
They can also be derived from the action of $\theta_2$ on the ramification locus in $A_2$
which, by piecing the local information together, leads to an invariant one-cycle in H$^1(\mathcal O^*)$,
i.e.~to the invertible sheaf $\mathcal O(R)$ on the quotient.)

Now, let $X$ be any supersingular K3 surface (or Kummer surface, for that thing) with Artin invariant 2. 
Then there is an isometry 
$$
\varphi: {\rm NS}({\rm GKm}(A)) \longrightarrow {\rm NS}(X)
$$
which maps  effective cycles to effective cycles (cf. Piatetskij-Shapiro--Shafarevich \cite{PSS},
Shioda \cite{Shi}).
Since we have 9 pairs of $(-2)$-curves $\{\psi_3(D'_i), \psi_3(F'_i)\}$ ($i = 1,  \ldots, 9$) 
in ${\rm NS}({\rm GKm}(A))$,
there are corresponding pairs $\{\varphi(\psi_3(D'_i)), \varphi(\psi_3(F'_i))\}$ ($i = 1,  \ldots, 9$)
in ${\rm NS}(X)$.
For the sake of simplicity, we set $D''_i = \varphi(\psi_3(D'_i))$ and 
$F''_i = \varphi(\psi_3(F'_i))$.
They are nonsingular rational curves which intersect transversally at a single point.
Reversing the above construction,
we blow up at each intersection point of $D''_i $ and $F''_i$.
Let 
$$
\psi'_{3}: \tilde{X} \longrightarrow X
$$ 
be the blowing-up. 
We denote by $G''_i$ ($i = 1,  \ldots, 9$) the exceptional curves while
we again use the same notation for the proper transforms of $D''_i$ and $F''_i$.
Then, $G''_i$, $D''_i$ and $F''_i$ form a triple 
with the same intersection numbers as before
(i.e.~$(G''_i, D''_i) = 1$,
$(G''_i, F''_i) = 1$, $(D''_i, F''_i) = 0$, $(D''_i)^2 = -3$, $(F''_i)^2 = -3$ 
and $(G''_i)^2 = -1$).
We set 
$$
R'' = \sum_{i = 1}^9 (D''_i  +2 F''_i).
$$
We can naturally extend the isometry $\varphi$ to 
$$
\tilde{\varphi} : {\rm NS}({\rm GKm}(A))\oplus_{i = 1}^9 {\bf Z}G'_i
\longrightarrow {\rm NS}(X)\oplus_{i = 1}^9 {\bf Z}G''_i,
$$
and we have $\tilde{\varphi}(R) = R''$.
Since $R$ is divisible by 3,  so is $R''$;
that is, there exists a divisor
$R'''$ on $\tilde{X}$ such that $R'' = 3R'''$.

\subsubsection{}

Using $R'''$, we can construct a cyclic covering 
$$
\pi' : \tilde{Y} \longrightarrow \tilde{X}
$$
of degree 3 with (smooth) branch locus the support of $R''$.
A local computation reveals that $\tilde Y$ can be taken to be smooth after a normalization
(which also follows from the general theory of triple covers in Miranda \cite{Miranda};
for K3 surfaces, see also Bertin \cite{Bertin}.)
By construction, $\tilde Y$ comes with an   induced  order 3 automorphism 
$\tilde\eta$ of $\tilde{Y}$ such that 
$$
\tilde{X} \cong \tilde{Y}/\langle \tilde\eta \rangle.
$$ 
Denoting pre-images by tildes, we have $\pi'^{-1}(G''_i) = \tilde{G''_i}$, $\pi'^{-1}(D''_i) = 3\tilde{D''_i}$ and 
$\pi'^{-1}(F''_i) = 3\tilde{F''_i}$ with $(\tilde{G''_i})^2 = -3$, $(\tilde{D''_i})^2 = -1$
and $(\tilde{F''_i})^2 = -1$. Contracting $\tilde{D''_i}$, $\tilde{F''_i}$, and subsequently $\tilde{G''_i}$,
we derive an algebraic surface $Y$, and we see that $\tilde\eta$ induces an automorphism $\eta$
of order 3 on $Y$. 

We will show that $Y$ is an abelian surface. 
The canonical divisor of $\tilde{X}$ is given by $K_{\tilde{X}} = \sum_{i= 1}^9G''_i$. 
Therefore, the canonical divisor of $\tilde{Y}$ is given by adding the ramification divisor:
$$
K_{\tilde{Y}} = \sum_{i = 1}^9\tilde{G''_i} + \sum_{i = 1}^92\tilde{D''_i} + \sum_{i = 1}^92\tilde{F''_i}.
$$
Hence, the canonical divisor of $Y$ is trivial. 

The Euler-Poincar\'e characteristic of $Y$ can be computed using a topological argument
for the involved coverings and blow-ups.
Essentially this works like over $\mathbf C$, except that we have to use \'etale cohomology with compact support.
For brevity, we omit the details leading to $\chi(Y)=0$.
%
Hence, considering the fact that $K_Y$ is trivial (and $p\geq 5$), we conclude that $Y$ is an abelian surface.

Since we have a dominant separable rational map from $Y$ to $X$ of degree 2,
by Lemma \ref{height} the height of the formal Brauer group $\Phi_Y$ 
is equal to the height of $\Phi_X$. Since $\Phi_X = \infty$, we have $\Phi_Y = \infty$.
Therefore, $Y$ is a supersingular abelian surface. By Lemma \ref{a-number of K3}
the a-number of $X$ is equal to 1. Therefore, by Lemma \ref{a-number}
the a-number of $Y$ is equal to 1, and $Y$ has Artin invariant 2.

\subsubsection{}

Recall that $Y/\langle \eta\rangle$ is
birational to $X$ and we have the following diagram:
$$
\begin{array}{ccccc}
&\tilde{Y} & \longrightarrow  & Y &\\
\swarrow& &   &   & \searrow    \\
\tilde Y/\langle\tilde\eta\rangle = \tilde{X}&  \longrightarrow & X &\longrightarrow   & Y/\langle \eta\rangle 
\end{array}
$$
From the diagram of exceptional curves
we see that the singularities of $Y/\langle \eta\rangle$ are of type $A_2$.
Since the a-number of $Y$ is equal to 1, the subgroup scheme $\alpha_p$ embeds uniquely into $Y$.
Therefore, $\eta$ preserves $\alpha_p$.
Let $P$ be a fixed point of $\eta$.
Let ${\mathcal O}_P$ be the local ring at the point $P$, and $m_P$ the maximal ideal.
We take the local parameter $s'$ in the direction of the subgroup scheme $\alpha_p$. Then, 
$\eta^*(s') = \gamma s' ~{\rm mod}~m_P^2$. Since $\eta$ is of order 3 and
$P$ is an isolated fixed point, we see that $\gamma$ is a primitive cube root of unity.
Setting 
$$
s = s' + \gamma^{-1}\eta^{*}(s') + \gamma^{-2}(\eta^{*})^2s',
$$ 
we see $\eta^*s = \gamma s$ and $s$ is  a nonzero element of $m_P/m_P^2$
by $p \neq 3$. Since the quotient singularity is of type $A_2$ and
the representation of ${\bf Z}/3{\bf Z}$ on $m_P/m_P^2$ is completely reducible,
we can take a element $t$ of $m_P$ such that 
$s, t$ form a basis of $m_P/m_P^2$ and such that $\eta^*s = \gamma s$,
$\eta^*t = \gamma^2 t$. 

Now, go to the quotient $Y' = Y/\alpha_{p}$. Then, $Y' \cong E \times E$ and 
$\eta$ induces an automorphism $\eta'$ of $Y'$. 
We may assume that $s^p, t$ give a local parameters
at a fixed point of $\eta'$. Therefore, the action of $\eta'$ at the fixed point
is given by $s^p \mapsto \gamma^2 s^p$,  $t \mapsto \gamma^2 t$.
Then, taking the Frobenius pull-back of these structures, we have the following
commutative diagram.
$$
\begin{array}{ccc}
       Y'^{(1/p)}    &\stackrel{\eta'^{(1/p)}}{\longrightarrow} &  Y'^{(1/p)} \\
              \downarrow  &        &  \downarrow  \\
          Y      &   \stackrel{\eta}{\longrightarrow}       &    Y \\
           \downarrow   &    & \downarrow  \\
         Y' & \stackrel{\eta'}{\longrightarrow} &   Y'.
\end{array}
$$
Here, all vertical arrows depict quotient morphisms by the subscheme $\alpha_p$ (suitably embedded),
and $Y'^{(1/p)} = Y' = E \times E$. We set $\tau = \eta''^{(1/p)}$. 
Let $Q$ be a fixed point of $\tau$, and
let ${\mathcal O}_Q$ be the local ring at the point $Q$, and $m_Q$ the maximal ideal.
Then, by our construction, we have local parameters $u, v$ of $m_Q$
such that the action of $\tau$ is given by $u \mapsto \gamma u$, $v \mapsto \gamma v$.
Since the cotangent space $m_Q/m_Q^2$ is isomorphic to 
${\rm H}^0(Y'^{(1/p)}, \Omega^1_{Y'^{(1/p)}})$, we see that the action of $\tau$
on the space ${\rm H}^0(Y'^{(1/p)}, \Omega^1_{Y'^{(1/p)}})$ is given 
by the multiplication by $\gamma$, which is a primitive cube root of unity.
Hence Lemma \ref{action} shows that
$Y'^{(1/p)}/\langle \tau \rangle$ is rational. 
Since the construction provides a purely inseparable morphism 
$$
Y'^{(1/p)}/\langle \tau \rangle \longrightarrow Y/\langle \eta \rangle
$$ 
of degree $p$,
and $Y/\langle \eta \rangle$ is birational  to $X$,
we conclude that $X$ is a Zariski surface as claimed.
\hfill\bsquare

\subsection{Proof of Theorem \ref{th} for $p\not\equiv 3$ mod $4$}
\label{ss}

If $p\equiv 3$ mod $4$, then the proof of Theorem \ref{th} proceeds very much along the lines of \ref{ss3},
except that there are a few subtleties to overcome since our automorphism does not have prime order
and Lemma \ref{action} does not apply.
As before, we start with a supersingular abelian surface $A$ with $\sigma=2$
endowed with the automorphism $\theta$ of order 4 from Section \ref{s:gen}.
Consider ${\rm GKm}(A)$, the minimal resolution of the quotient $A/\langle \theta\rangle$
with singularities of types $4A_3+6A_1$.
Thus GKm$(A)$ carries natural configurations of smooth rational curves
$C_i, D_i, E_i \; (i=1,2,3,4)$, forming $A_3$ root lattices, and 6 disjoint $(-2)$-curves $F_i \; (i=1,\dots,6)$.
In particular, GKm$(A)$ contains an effective 4-divisible divisor 
\[
R = \sum_{i=1}^4 (C_i+2D_i+3E_i) + 2\sum_{i=1}^6 F_i,
\]
but the corresponding cover $A_0 \to$ GKm$(A)$ is not smooth since supp$(R)$ is not.
We will see momentarily how to overcome this without any additional blow-ups
(one of the advantages over the direct approach from \ref{ss3}).

Now consider a supersingular K3 surface $X$ of Artin invariant 2.
Then, as before, NS$(X)$ contains the same configuration
of $(-2)$-curves, and the same 4-divisible divisor $R$ (using the same notation as for GKm$(A)$).
We search for a smooth birational model of the corresponding cyclic degree 4 cover $Y_0$.
To this end, we first consider the smooth degree 2 cover 
\[
\tilde W \to X
\]
corresponding to the  2-divisible branch divisor $\sum_{i=1}^4(C_i+E_i)$ with smooth support.
Clearly the $C_i, E_i$ pull-back to $(-1)$-curves $C_i', E_i'$ on $\tilde W$
while 
\[
K_{\tilde W} = \sum_{i=1}^4(C_i'+E_i').
\]
Contracting these disjoint $(-1)$-curves, we obtain a smooth surface $W$ with
\[
K_W = 0, \;\;\; \chi(W)=24,
\]
thus a K3 surface. This comes equipped with $(-2)$-curves $D_i'$ mapping to $D_i$,
and with effective $(-4)$-divisors mapping $2:1$ to the $F_i$.  
It follows that these decompose into two disjoint $(-2)$-curves $F_{i,1}+F_{i,2}$ each.

On $\tilde W$, $R$ pulls back to twice the divisor $R'+2\sum_{i=1}^9 E_i'$ where
\[
R' = \sum_{i=1}^4 (C_i' + D_i' + E_i') + \sum_{i=1}^6 (F_{i,1} + F_{i,2})
\]
is still 2-divisible  by construction; the corresponding cover remains birational to $Y_0$.
Push forward of $R'$ to the K3 surface $W$ yields the 2-divisible divisor
\[
R'' = \sum_{i=1}^4 D_i'  + \sum_{i=1}^6 (F_{i,1} + F_{i,2})
\]
which again is smooth. As before, the corresponding cover is the blow-up $\tilde Y$ of an abelian surface $Y$,
this time in 16 points,
and indeed $W$ is the Kummer surface of $Y$.
(To rule out (quasi-)bielliptic surfaces (in characteristic $3$), one may note that
by construction, $W$ is supersingular and hence $\rho(Y)\geq 6$.)
The whole construction endows $Y$ with an automorphism $\eta$ of order 4 whose quotient is birational to $X$.
We sketch the resulting maps in the following diagram:
$$
\begin{array}{ccccc}
&& \tilde Y & \to & Y\\
&& \downarrow && \downarrow \\
\tilde W & \to & W & \to & Y/\langle\eta^2\rangle \\
\downarrow &&&& \downarrow\\
X && \to && Y/\langle\eta\rangle\\
\end{array}
$$

Now we can proceed exactly as before, with intermediate step from $X$ to $Y$ going through $W$,
to deduce that all these varieties are supersingular with a-number 1.
It follows that $\alpha_p$ admits a unique embedding into $Y$
which is thus compatible with the action of $\eta$.
This induces an order 4 automorphism $\eta'$ on the quotient $Y'=Y/\alpha_p$.
However, we cannot infer from Lemma \ref{action} hat $Y'/\langle\eta'\rangle$ is rational,
so we have to pursue a different line of reasoning.
To this end, we once again factorize the quotient map.
Namely, we first consider the quotient 
$$
V'=Y'/\langle\eta'^2\rangle
$$
Note that $\eta^2=-$id on $Y$, so the same holds for $\eta'^2$ on $Y'$.
That is, $V'$ is birational to Km$(Y')$, a K3 surface.
Then $\eta'$ induces an involution $\imath$ on $V'$.
Since $\imath$ kills the regular 2-form on $V'$ by construction,
the quotient 
$$
V=V'/\langle\imath\rangle\sim Y'/\langle\eta'\rangle
$$ 
cannot be (birationally) K3 again.
It follows that $V$ is either rational or Enriques.
But then $V$ admits a purely inseparable map of degree $p>2$ to the K3 surface $X$ by construction,
so $V$ cannot have fundamental group ${\mathbf Z}/2{\mathbf Z}$.
Hence $V$ is rational, and $X$ is Zariski as claimed.
This completes the proof of Theorem \ref{th}.
\hfill\bsquare

\begin{remark}
It may be feasible to pursue an alternative approach to prove Theorem \ref{th}
based on the results of Blass--Levine \cite{BL},
choosing a suitable polarization etc.
However, given the explicit geometric arguments
which supersingular Kummer surfaces lend themselves to,
we decided to stick to the above reasoning.
\end{remark}

\section{Supersingular K3 surfaces}

In order to deal with other supersingular K3 surfaces
and compute their Artin invariants,
we need a little preparation on the lattice theoretic side.
To this end, we assume that $S$ is a supersingular K3 surface,
endowed with a certain sublattice $L$ embedding into the N\'eron-Severi group ${\rm NS}(S)$.
For instance, $S$ could be the minimal resolution of some singular surface
with $L$ generated by the exceptional curves above the singularities.


\begin{th}\label{Artin invariant} 
In the notation above, let $n$ be the rank of $L$.
We asssume that the discriminant
of $L$ is prime to the characteristic $p$. Then the Artin invariant
$\sigma$ of $S$ is smaller than or equal to $(22 - n)/2$.
\end{th}

\proof{We denote by ${\rm NS}(S)^{*}$ the dual lattice of ${\rm NS}(S)$, and likewise for $L$ etc.
Then, by the result by M. Artin (\cite{A}), ${\rm NS}(S)^{*}/{\rm NS}(S)$ is 
a $p$-elementary group, and $\mid {\rm NS}(S)^{*}/{\rm NS}(S) \mid = p^{{2\sigma}}$.
Without loss of generality, we assume that 
$L$ embeds primitively into ${\rm NS}(S)$,
for else we could continue with the primitive  closure of $L$ inside ${\rm NS}(S)$
which will still have discriminant prime to $p$.
Note that by definition, its orthogonal complement $L^\perp$ embeds primitively  into ${\rm NS}(S)$ as well.

Consider the embedding
\[
L\oplus L^\perp \hookrightarrow {\rm NS}(S)
\]
of finite index $m$, say.
The way how $L$ and $L^\perp$ glue together 
is encoded in an isomorphism of subgroups of the discriminant groups,
\[
L^*/L \supseteq H_1 \cong H_2 \subseteq (L^\perp)^*/L^\perp
\]
such that the induced intersection form (modulo $2\mathbf Z$) agrees up to sign.
We use two related properties:
on the one hand, from lattice theory,
\[
|H_1| = |H_2| = m;
\]
on the other hand, as a subgroup,
\[
|H_1|  \mid |L^*/L| = |{\rm disc} L|.
\]
In particular, our assumption implies that $m$ is prime to $p$.
From this we aim to deduce that (the $p$-part of) ${\rm NS}(S)^{*}/{\rm NS}(S)$
is fully captured in $ (L^\perp)^*/L^\perp$, with structure unchanged.
To this end, 
consider the sequence of ${\bf Z}$-modules:
\begin{eqnarray*}\label{sequence}
  L^{\perp}\subset L\oplus L^{\perp}\subset {\rm NS(S)} \subset {\rm NS(S)}^{*}\subset (L\oplus L^{\perp})^{*}
  \cong L^{*}\oplus (L^{\perp})^{*} \supset (L^{\perp})^{*}
\end{eqnarray*}
where the finite index inclusions in the middle have index $m, p^{2\sigma}$ and $m$, respectively.
It follows that $m^{2}p(L\oplus L^{\perp})^{*}\subset L\oplus L^{\perp}$,
and thus, restricting to the orthogonal summand $L^\perp$, also $m^{2}p(L^{\perp})^{*} \subset L^{\perp}$. Since ${\rm rank}~L^{\perp} = 22 -n$, 
we infer from the elementary divisor theorem that
\begin{equation}\label{divisor}
 p^{23 - n}\not| \;{\rm disc} \; L^\perp.
\end{equation}
Putting everything together, we use
\[
-p^{2\sigma} = {\rm disc \; NS}(S) = \frac{({\rm disc} \; L)({\rm disc} \; L^\perp)}{m^2} 
\]
to deduce, from our assumption that ${\rm disc} \; L$ and thus $m$ is prime to $p$,  and from (\ref{divisor}), that
$2\sigma \leq 22 -n$ as claimed.
}

\begin{ex}\label{Kummer}
Let $k$ be an algebraically closed field of characteristic $p\neq 2, 3$,
and $A$ an supersingular abelian surface over $k$. Then $A$ has a principal polarization $\Theta$.
Since $p \neq 2, 3$, we can choose a nonsingular curve of genus 2 as the principally polarization
(cf. Ibukiyama-Katsura-Oort \cite{IKO}, and Ogus \cite{O}). Consider the linear system $\mid 2\Theta \mid$.
Then the associated rational map $\varphi_{\mid 2\Theta \mid}$ is a morphism (cf. Mumford \cite{M}),
and, as is well-known, the image of $\varphi_{\mid 2\Theta \mid}$
is a quartic surface with 16 rational double points of type $A_{1}$ in the projective plane ${\bf P}^{3}$, 
which is isomorphic to the quotient surface $A/\langle \iota \rangle$. Here, $\iota$ is 
the inversion of $A$.
The minimal resolution of the surface $A/\langle \iota \rangle$
is the Kummer surface ${\rm Km}(A)$. We denote by 
$$
\pi : {\rm Km}(A) \longrightarrow A/\langle \iota \rangle
$$ 
the resolution.
We take a generic hyperplane section $H$ of $A/\langle \iota \rangle$
and pull-back $D = \pi^{*}H$. Then, $D$ does not intersect
the exceptional divisors and $D^{2}= 4$, which is prime to $p$. 
Consider the lattice $L\subset$ NS$({\rm Km}(A))$ 
generated by the exceptional curves and $D$.
Then, we have {\rm rank} $L$ = 17. By the same argument as in Theorem \ref{Artin invariant},
we see that the Artin invariant satisfies
$$
\sigma ({\rm Km} (A))\leq \frac{22 -17}{2}.
$$
 Thus we obtain an alternative reasoning for the well-known result 
$\sigma ({\rm Km} (A)) \leq 2$ (cf. Ogus \cite{O}).
\end{ex}

\begin{ex}
Analogous arguments apply to the generalized Kummer surfaces from Section \ref{s:gen}
to prove that they have Artin invariant $\sigma\leq 2$.
\end{ex}

We continue with another application to Kummer surfaces
which is a kind of converse of Example \ref{Kummer}.
Recall that
Nikulin showed that a complex K\"ahler K3 surface $X$ is a Kummer surface if and only if there exist 16 nonsingular rational
curves on X which do not intersect each other (cf. Nikulin \cite{N}). For supersingular K3 surfaces,
as an application to Theorem \ref{Artin invariant}, we have the following:

\begin{th} Let $X$ be a supersingular K3 surface with a divisor D such that $D^2$ is prime to the characteristic $p$.
Assume that there exist 16 nonsingular rational curves $E_i$ which do not intersect each other, and assume
that $(D, E_i) = 0$ $(i = 1, 2, \cdots, 16)$. Then, $X$ is a Kummer surface.
\end{th}

\proof{We consider the lattice $L\subset{\rm NS}(X)$ generated by $D$ and $E_i$'s. The lattice $L$ is of rank 17 and
the discriminant is prime to $p$. Therefore, by Theorem \ref{Artin invariant}, the Artin invariant $\sigma$
is smaller than or equal to $(22 -17)/2$. Therefore, we have $\sigma \leq 2$. Hence, $X$ is
a Kummer surface (cf. Ogus \cite{O}, Shioda \cite{Shi}).
}

We proceed by giving a direct construction covering the remaining part of Theorem \ref{thm} \textit{(ii)}.
We emphasize that this does not require any further machinery;
in the next section it will be generalized in the context of elliptic surfaces
(Lemma \ref{lem:ell} etc).

\begin{lm}\label{ell=5}
Assume $\ell = 5$ in $(\ref{equation1})$. Then, $(C \times C)/\langle \tau \times \tau^2\rangle$ 
is birational to a K3 surface $S$, and
$(C \times C)/\langle \tau \times \tau^3\rangle$ is  birational 
to the same K3 surface.
\end{lm}
\proof{The group $G = \langle \tau \times \tau^{2} \rangle$ acts on
the function field of $C \times C$ via its natural action on $k(x, y, x_1, y_1)$. We set
$$    
     z = xx_{1}^2.
$$
Then, $z$ is invariant under $G$ and the invariant field $k(C \times C)^G$
 is given by $k(y_{1}, y, z)$ with the equation 
$z^5 = (y^2 + y) (y_1^2 + y_1)^2$.
We set
$$
    w = y(y_1^2 + y_1)
$$
Then the relation translates as
\begin{equation}\label{K3} 
     w^2 + w(y_1^2 + y_1) = z^5
\end{equation}
which gives a birational equation for the quotient surface $S$.
Outside characteristic $2$, $S$ is thus birational to the double cover of ${\bf P}^2$
branched along the sextic curve
\begin{equation}\label{ramification}
C:\;\;\; z^{5}u + \frac{u^2(y_1^2+y_1u)^2}{4} = 0
\end{equation}
where $y_1,z,u$ denote homogeneous coordinates of ${\bf P}^2$.
Obviously $C$ is reducible, but the singularities are only isolated rational double points.
It follows that the minimal resolution of the double cover is a K3 surface as claimed.

In characteristic $2$, (\ref{K3}) still defines a separable double cover of ${\bf P}^2$,
but due to the presence of wild ramification,
the branch locus degenerates to the cubic curve $(y_1^2 + y_1u)u$.
Yet the singularities and the underlying invariants are preserved,
so we obtain a K3 surface as before.

Since $\tau$ is of order $5$, $(\tau \times \tau^2)^3 = (\tau^3 \times \tau)$ is 
a generator of the group $\langle \tau \times \tau^2\rangle$.
Therefore, by exchanging the components of $S$, we have 
an isomorphism from $(S)/\langle \tau \times \tau^2\rangle$
to $(S)/\langle \tau \times \tau^3\rangle$.
This concludes the proof of Lemma \ref{ell=5}.
}

We are now in the position to prove the remaining part of Theorem \ref{thm} \textit{(ii)}.

\begin{th}
\label{thm:5}
Assume $p \equiv 2 ~{\rm or}~ 3~ ({\rm mod}~5)$.
Then the K3 surface $S$ from Lemma \ref{ell=5} is Zariski 
with Artin invariant 2. 
\end{th}

\proof{By Proposition \ref{Zariski} and Lemma \ref{ell=5}, $S$ is a Zariski K3 surface.
In order to exhibit a suitable sublattice $L$ of NS$(S)$,
we study the singularities of the double covering of ${\bf P}^2$ from the proof of Lemma \ref{ell=5}.
In the affine chart (\ref{K3}), there are two singularities at $(w,y_1,z)=(0,0,0), (0,-1,0)$;
visibly, each is a rational double point of type $A_4$.
The chart $y_1\neq 0$ with affine equation
\[
 w^2 + w(1 + u)u = uz^5
 \]
has another rational double point at $(w,u,z)=(0,0,0)$, this time of type $A_9$.
The minimal resolution $S$ is thus endowed with the sublattice $L\subset{\rm NS}(S)$
generated by the exceptional curves above the singularities.
Presently $L$ has rank $17$ and discriminant $250$.
In particular, the discriminant is prime to $p$ if $p>2$,
so using Theorem \ref{Artin invariant},
we see that the Artin invariant $\sigma$ of $S$ is smaller than or equal to 2.
To establish the same claim in characteristic two,
it suffices to enhance the lattice $L$ by the classes of the strict transforms of the two 'lines' $\{w=z=0\}$
and $\{w=u=0\}$.
One easily checks that the resulting overlattice $L'$ has rank $18$ and discriminant $-5$,
so we conclude $\sigma\leq 2$ as before.

In order to prove the equality $\sigma=2$,
we appeal to work of Jang \cite{J} studying the non-symplectic index $N$ of supersingular K3 surfaces.
This is defined as the size of the image of the natural representation
\[
{\rm Aut}(S) \to {\rm GL}(H^0(S,\Omega_S^2)).
\]
In detail, Jang proves in characteristic $p>3$ that
\begin{itemize}
\item
the supersingular K3 surface with $\sigma=1$ has $N=p+1$, and
\item
the supersingular K3 surfaces with $\sigma=2$ have $N=2$ or, at a unique moduli point, $N=p^2+1$.
\end{itemize}
Presently, the automorphism $(w,z,y_1)\mapsto (w,\zeta z,y_1)$ induced by $\tau$
acts primitively of order $5$ on the regular two-form $dz\wedge dy_1/w$,
so the initial congruence assumption for $p$  implies $N=p^2+1$ and $\sigma=2$ as claimed.

In characteristics $2$ and $3$ where Jang's results are not valid,
one can complete the proof without difficulty using elliptic fibrations.
We postpone their treatment until Remark \ref{rem:2+3}.
}

Note that the proof of Theorem \ref{thm} \textit{(ii)} is now complete
since the two characteristics so far missing from the proof of Theorem \ref{thm:5}
are covered by Theorem \ref{2}.

%

\section{Zariski elliptic surfaces}

In this section, we continue to argue with the curve $C$ from (\ref{equation1})
with automorphism $\tau$ for odd $\ell$
and consider the quotient
\[
S = (C\times C)/\langle\tau\times\tau^i\rangle \;\;\; (i\in\{1,\dots,\ell-1\}).
\]
For starters, we briefly leave the restricted area of K3 surfaces:

\begin{lm}
\label{lem:ell}
Let $i=(\ell-1)/2$. If $p\equiv i$ mod $\ell$ or $p\equiv i^{-1}$ mod $\ell$,
then $S$ is a Zariski surface admitting an elliptic fibration. 
\end{lm}

\proof{
It is immediate that $S$ is still birationally given by the degree $\ell$ analogue of the affine equation (\ref{K3}):
\begin{eqnarray}
\label{eq:S}
S: \;\;\; w^2 + w(y_1^2+y_1) = z^\ell.
\end{eqnarray}
Interpreting this as a cubic over $k(z)$, we obtain the claim
(with sections at $\infty$, indeed).
}

\begin{remark}
Of course, there are Zariski elliptic surfaces in the literature,
but to our knowledge mostly arising by purely inseparable base change from a rational elliptic surface
(Shioda \cite[Ex.~4.2]{Shioda-some}, Katsura \cite{K0}) -- just like in Example \ref{ex:13}. 
\end{remark}

Transferring (\ref{eq:S}) to a Weierstrass form is easily achieved:
homogenize as a cubic in ${\bf P}^2_{k(z)}$ by a variable $u$, say
and switch to the affine chart $w\neq 0$ to derive
\[
S: \;\;\; y_1^2+uy_1 = z^\ell u^3
-u.
\]
Then multiplying the equation by $z^{2\ell}$ and dividing variables by $z^\ell$,
we arrive at the normalized Weierstrass form
\begin{eqnarray}
\label{eq:S2}
S: \;\;\; y_1^2+uy_1 =  u^3
-z^\ell u.
\end{eqnarray}
Note the two-torsion section at $(0,0)$, and the automorphism $\tau$ induced by the action $z\mapsto \zeta z$ on the base.

\begin{lm}
\label{lem:5+7}
In the above setting,
$S$ is birational to a K3 surface if and only if $\ell=5$ or $7$.
\end{lm}

\proof{
This is a standard argument using the theory of elliptic fibrations,
see e.g.~\cite{SSh}.
To give some details, we compute the basic invariants of (the Kodaira--N\'eron model of) $S$,
starting from the zero irregularity (which follows from the base curve being ${\bf P}^1$).
The fibration (\ref{eq:S2}) has singular fibers of Kodaira types
\[
I_{2\ell}/z=0, \;\;\; I_1/64z^\ell=-1,\;\;\; 
~
\left\{
{\begin{array}{cc}
III/\infty, & \ell \equiv 3\; {\rm mod} \; 4,\\
III^*/\infty, & \ell \equiv 1\; {\rm mod} \; 4,
\end{array}}
\right.
\]
except that in characteristic two, the $I_1$ fibers are absorbed by the wild ramification at $\infty$.
It follows that $S$ has Euler-Poincar\'e characteristic $e(S)=3\ell+3$ resp.~$e(S)=3\ell+9$ and
geometric genus $\lfloor\ell/4\rfloor$.
Hence $S$ is a K3 surface exactly for $\ell=5$ and $7$ as claimed.
}

\begin{remark}
\label{rem:2+3}
As a first application, we explain how to infer that for $\ell=5$ and $p=2$ or $3$,
the Artin invariant of $S$ cannot be $\sigma=1$ (as stated in Theorem \ref{thm:5}).
To see this, it suffices to go through the classification of elliptic fibrations on the superspecial K3 surface
in each characteristic:
by inspection neither \cite{ES} nor \cite{Sengupta} lists a fibration with both given reducible fibers of type $III^*$ and $I_{10}$.
\end{remark}

We now proceed to prove Theorem \ref{thm} \textit{(iii)}.

\begin{th}
\label{thm:7}
Let $\ell=7$ and $p\equiv 3,5$ mod $7$.
Then $S$ is a Zariski K3 surface with Artin invariant $\sigma=3$.
\end{th}

\proof{
The surface $S$ is Zariski by Theorem \ref{Zariski}
and K3 by Lemma \ref{lem:5+7}.
We proceed by exhibiting a suitable sublattice $L\subset$ NS$(S)$.
To this end, consider the lattice $L'$ generated by fiber components and zero section,
\[
L' = U \oplus A_1 \oplus A_{13},
\]
where $U$ denotes the hyperbolic plane generated by zero section $O$ and general fiber.
The two-torsion section provides an index $2$ overlattice $L'\subset L\subset$ NS$(S)$
of rank $16$ and discriminant $-7$.
Hence Theorem \ref{Artin invariant} applies to show that $\sigma\leq 3$.
If $p>3$, then we conclude using the non-symplectic index following Jang as before.
For $p=3$, in contrast, we pursue a direct approach by exhibiting a full set of generators of NS$(S)$
using the theory of Mordell--Weil lattices after Shioda \cite{ShMW}.
In practice, we search for a section $P$ of small height.
This soon leads to $P$ being integral (i.e.~disjoint from the zero section $O$)
and meeting both reducible fibers in components adjacent to the identity component.
This means that $P=(tU,tV)$ for polynomials $U,V\in k[z]$ of degree $2$ resp.~$4$
with $t\not| U$.
In fact, the special shape of the Weierstrass form (\ref{eq:S2}),
in particular the presence of the two-torsion section,
implies that $U$ has to be a square in $k[z]$.
Given this, one can solve directly for $P$ to find, uniquely up to  symmetry,
\[
P = (-t(t+1)^2, t^2(t+1)(t-1)^2).
\]
Comparing the seven sections $P_j = \tau^j P \; (j=0,\dots,6)$, we find that any two of them intersect
transversally in a single point. Hence the height pairing evaluates as
\begin{eqnarray*}
h(P_j) & = & 4 + 2 \underbrace{(P_j.O)}_{=0}  - \frac {13}{14}-\frac 12\; \;\; = \;\frac{18}7 \;\;\; (0\leq j\leq 6)\\
\langle P_j, P_m\rangle  & = & 2 -  \underbrace{(P_j.P_m)}_{=1}  - \frac {13}{14}-\frac 12 \;  = \; -\frac 37 \;\; \;(0\leq j\neq m\leq 6)
\end{eqnarray*}
where the correction terms are read off from the fiber components met.
From the resulting Gram matrix, we infer that the $P_i$ generate a sublattice $M$ of the Mordell--Weil lattice of $S$
of rank $6$ and discriminant $3^6/7$,
in perfect agreement with the fact $\sigma\leq 3$.
Proving equality thus amounts to showing that $M$ equals the full Mordell--Weil lattice,
i.e.~that there cannot be any divisibilities among the given sections.
This can be verified in multiple ways,
for instance using the fact that the automorphism $\tau$ makes $M$  an (irreducible) ${\bf Z}[\zeta]$-module of rank one.
Hence a single divisibility would cause several independent others -- too many in fact for the discriminant of NS$(S)$ to stay integral.
}

\begin{remark}
One could also argue without appealing to Mordell-Weil lattices,
just using intersection numbers and the rank formula often attributed to Shioda-Tate \cite[Cor.~1.5]{Sh-EMS}.
The above reasoning, in contrast, seems more streamlined and conceptual.
\end{remark}

\begin{remark}
\label{rem:ell}
Similar results can be derived for the  Zariski elliptic surfaces with $e>24$ from Lemma \ref{lem:ell}.
In particular, the length of the discriminant group (generalizing twice the Artin invariant) is always bounded by $\ell-1$
(and even smaller when $\ell$ is not prime).
We emphasize that our approach is not limited to characteristics satisfying the standard condition
\begin{eqnarray}
\label{eq:nu}
\exists\, \nu: \;\; p^\nu \equiv -1 \;\;{\rm mod} \; \ell
\end{eqnarray}
from the Fermat surface case \cite{Shioda-p=3}, \cite{KS}.
Indeed, for $\ell=11$, for instance, we obtain (non-rational) Zariski surfaces
in characteristics congruent to $5,9$ modulo $11$ which do not satisfy (\ref{eq:nu}).
\end{remark}

As a  supplement to Theorem \ref{thm} \textit{(i)}, 
we conclude this paper by providing a Zariski K3 surface in characteristic $p=13$ of Artin invariant $\sigma=1$
in the vein of Shioda \cite[Ex.~4.2]{Shioda-some}.

\begin{ex}
\label{ex:13}
Assume $p = 13$.
Consider the rational elliptic surface given in Weierstrass form
\[
 y^2 = x^3 + tx - t.
\]
It has singular fibers of Kodaira type $III^*$ at $\infty$,
$II$ at $t=0$ and $I_1$ at $t=-27/4$.
Applying the purely inseparable base $t=s^{13}$,
we obtain an elliptic  K3 surface $X$ with the same additive fibers,
but $I_1$ replaced by $I_{13}$.
By construction, $X$ is Zariski and furnished with a sublattice
\[
L = U \oplus A_{12} \oplus E_7\subset {\rm NS}(X)
\]
of rank $21$ and discriminant $26$.
While this is not relatively prime to the characteristic,
Theorem \ref{Artin invariant} is easily adjusted to prove
that $X$ has Artin invariant $\sigma\leq 1$.
Hence equality holds.
(Alternatively this can be inferred from the section $P=(1,1)$ of height $1/2$ on the rational elliptic surface.)
%
\end{ex}

\noindent
{\bf Acknowledgement.}
We thank T.\ Shioda for helpful comments and discussions.

\vspace{0.5cm}
\noindent
T.\ Katsura: Faculty of Science and Engineering, Hosei University,
Koganei-shi, Tokyo 184-8584, Japan

\noindent
E-mail address: toshiyuki.katsura.tk@hosei.ac.jp

\smallskip

\noindent
M.\ Sch\"utt: 
Institut f\"ur Algebraische Geometrie, Leibniz Universit\"at
  Hannover, Welfengarten 1, 30167 Hannover, Germany, and
  Riemann Center for Geometry and Physics, 
  Appelstrasse 2, 30167 Hannover, Germany

\noindent
E-mail address: 
schuett@math.uni-hannover.de

\end{document}